\documentclass[11pt,leqno]{amsart}

\usepackage{amssymb}
\usepackage{amsmath}
\usepackage{amsthm}
\usepackage{latexsym}
\usepackage{amscd, xypic}
\usepackage[dvipsnames]{xcolor}
\usepackage{tikz-cd} 
\tikzcdset{arrow style=tikz}
\usepackage{enumitem}
\setlist[itemize]{leftmargin=*}
\setlist[enumerate]{leftmargin=*}

\begin{document}

\newtheorem{thm}{Theorem}
\newtheorem{prop}[thm]{Proposition}
\newtheorem{lem}[thm]{Lemma}
\newtheorem{cor}[thm]{Corollary}
\newtheorem*{THM}{Theorem}
\newtheorem*{LEM}{Lemma}
\newtheorem*{PRO}{Proposition}
\theoremstyle{definition}
\newtheorem{defi}[thm]{Definition}
\newtheorem{ex}[thm]{Example}
\newtheorem{rem}[thm]{Remark}

\newcommand{\ob}{\mbox{ob}}
\newcommand{\mor}{\mbox{mor}}
\newcommand{\iso}{\mbox{Iso}}
\newcommand{\id}{\mbox{id}}
\newcommand{\G}{\mathcal {G}}
\newcommand{\R}{\mathbb{R}}
\newcommand{\N}{\mathbb{N}}
\newcommand{\si}{\sigma}
\newcommand{\rh}{\rho}
\newcommand{\ta}{\tau}
\newcommand{\rmod}{$R$\mbox{-mod}}
\newcommand{\grmod}{\mbox{$G$-$R$-mod}}
\newcommand{\HOM}{{\rm HOM}}
\newcommand{\Hom}{\mbox{Hom}}
\newcommand{\END}{{\rm END}}
\newcommand{\Ker}{\mbox{Ker}}
\newcommand{\Image}{\mbox{Im}}
\newcommand{\m}{^{-1}}

\newcommand{\GR}{\mbox{$\G$-$R$}}
\newcommand{\AbG}{\mbox{$\mbox{Ab}_{\G}$}}
\newcommand{\Ab}{\mbox{Ab}}
\newcommand{\rmd}{\mbox{$R$-md}}
\newcommand{\mdr}{\mbox{md-$R$}}
\newcommand{\rmds}{\mbox{$R$-md-$S$}}
\newcommand{\smod}{\mbox{$S$-mod}}
\newcommand{\rmods}{\mbox{$R$-mod-$S$}}
\newcommand{\modr}{\mbox{mod-$R$}}
\newcommand{\gsmod}{\mbox{$\G$-$S$-mod}}
\newcommand{\gmodr}{\mbox{$\G$-mod-$R$}}
\newcommand{\grmods}{\mbox{$\G$-$R$-mod-$S$}}
\newcommand{\grmodr}{\mbox{$\G$-$R$-mod-$R$}}
\newcommand{\gmods}{\mbox{$\G$-mod-$S$}}
\newcommand{\grmd}{\mbox{$\G$-$R$-md}}
\newcommand{\gsmd}{\mbox{$\G$-$S$-md}}
\newcommand{\gmdr}{\mbox{$\G$-md-$R$}}
\newcommand{\gmds}{\mbox{$\G$-md-$S$}}
\newcommand{\grmds}{\mbox{$\G$-$R$-md-$S$}}

\newcommand{\rumod}{R\text{-\textup{\textbf{umod}}}}
\newcommand{\sumod}{S\text{-\textup{\textbf{umod}}}}
\newcommand{\modgr}{\mbox{{\bf mod}-}\G\mbox{-$R$}}
\newcommand{\grumod}{\G\textup{-}\rumod}
\newcommand{\gsumod}{\G\textup{-}\sumod}
\newcommand{\umodgr}{\G\text{-\textup{\textbf{umod}}}\textup{-}R}
\newcommand{\umodgs}{\G\text{-\textup{\textbf{umod}}}\textup{-}S}

\title[Primitives of continuous functions 
via polynomials]{Primitives of continuous 
functions \\ via polynomials}

\author[P. Lundstr\"{o}m]{Patrik Lundstr\"{o}m}
\address{University West,
Department of Engineering Science, 
SE-46186 Trollh\"{a}ttan, Sweden}
\email{{\scriptsize patrik.lundstrom@hv.se}}



\begin{abstract}
We present an elementary self-contained
folkloristic proof, using limits 
of primitives of Bernstein polynomials,
for the existence of primitive functions
of continuous functions defined on the 
unit interval.
\end{abstract}

\maketitle

In most standard calculus books the existence
of primitives of continuous functions is proved 
using limits of Riemann sums 
(see e.g. \cite{adams2006,stewart2015}).
Therefore, to most students, and even to some
instructors, it might seem like this is the 
only way to show this.
The aim of this note is to highlight an elementary 
but completely different proof of this fact
which combines on the one hand the so called 
Stone-Weierstrass theorem 
\cite{perez2008,stone1948}
on uniform approximation 
\linebreak
of continuous functions
by polynomials, and
on the other hand a classical result from calculus
on the existence of limits of differentiated 
sequences of 
\linebreak
functions \cite[Thm. 9.13]{apostol1974}.
Somewhat simplified, the idea of the proof is to
\linebreak 
show that primitives 
to a sequence of polynomials converging to a given
\linebreak 
continuous function, converges to a 
primitive function of the given function.

Throughout the article, we use the following notation.
Let $I$ and $J$ denote, respectively, 
the closed unit interval $[0,1]$ and the open unit 
interval $(0,1)$.
All functions are supposed to be real-valued.
For a bounded function $g$ on $I$
we put $|| g || = {\rm sup} \{ |g(x)| \}_{x \in I}$.
For future reference, we recall the following:

\begin{PRO}
Suppose that a sequence $g_1,g_2,g_3,\ldots$ of 
continuous functions on $I$ is 
uniformly convergent, that is
$|| g_m - g_n || \to 0$ 
as $\mbox{min} (m,n) \to \infty$.
\linebreak
Then the limiting function
$g(x) := \lim_{n \to \infty} g_n(x)$ 
is continuous on $I$.
\end{PRO}

\noindent
\emph{Proof of proposition.}
Take $x$ in $I$ and $\epsilon > 0$. 
Choose an integer $n \geq 1$ so that
$|| f_n - f || < \epsilon/3$.
Since $f_n$ is continuous on $I$ there is 
$\delta > 0$ such that $| f_n(x) - f_n(y) | < \epsilon/3$
whenever $y$ belongs to $I$ 
and $|x - y| < \delta$.
Therefore,
\linebreak
for such a real number $y$, 
the triangle inequality now implies that 
\linebreak
$| f(x) - f(y) | \leq | f(x) - f_n(x) | +
| f_n(x) - f_n(y) | + | f_n(y) - f(y) | <
3 \cdot \frac{\epsilon}{3} = \epsilon$.
\qed

\vspace{2mm}

For the rest of the article, $f$ denotes a fixed
continuous function defined on $I$.
Let $m$ and $n$ be integers. Define
$\smash{p_{m,n}(x) = {n \choose m} x^m (1-x)^{n-m}}$
where we use the standard convention that 
$\smash{{n \choose m} = 0}$ whenever $m > n$ or $m < 0$. 
\linebreak
If $n \geq 1$, then we put 
$\smash{F_n(x) = \sum_m f(m/n)/(n+1)  
\sum_{j=m+1}^{n+1} p_{j,n+1}(x)}$.
\linebreak
We wish to show the following:

\begin{THM}
The sequence of polynomials $F_1,F_2,F_3,\ldots$
converge uniformly on $I$ to a 
continuous function $F$ which is differentiable
on $J$ with $F' = f$. 
\end{THM}

To this end, we first prove the following:

\begin{LEM}
Let $m,n$ be integers with $n \geq 1$.
The following assertions hold:
\begin{enumerate}

\item $p_{m,n}'(x) = 
n \big( p_{m-1,n-1}(x) - p_{m,n-1}(x) \big)$

\item $1 = \sum_{m=0}^n p_{m,n}(x)$

\item $nx = \sum_{m=0}^n m p_{m,n}(x)$

\item $nx(1-x) = \sum_{m=0}^n (nx-m)^2 p_{m,n}(x)$

\item $F_n(0) = 0$ and $F_n' = f_n$ where  
$f_n(x) = \sum_{m=0}^n f 
\left( m/n \right) p_{m,n}(x)$

\item The sequence of polynomials $f_1,f_2,f_2,\ldots$
converge uniformly to $f$ on $I$. 

\end{enumerate}
Note that the polynomial $f_n$ above is called the 
$n$th Bernstein polynomial \cite{bernstein1912}.
\end{LEM}

\noindent
\emph{Proof of lemma.}
The equalities (2)-(4) above can be shown quite easily 
by using 
the binomial theorem (see \cite[p. 7]{burkill1959}).
Here we give a direct argument.
\begin{enumerate}

\item By the product rule for derivatives it follows that
{\footnotesize
\[
p_{m,n}'(x) = {n \choose m} m x^{m-1}(1-x)^{n-m} - 
{n \choose m} (n-m) x^m (1-x)^{n-m-1} 
\]
\[
=
n {n-1 \choose m-1} x^{m-1}(1-x)^{n-m} -
n {n-1 \choose m} x^m (1-x)^{n-m-1} =
n \big( p_{m-1,n-1}(x) - p_{m,n-1}(x) \big).
\]
}

\item Put $\alpha(x) = \sum_{m=0}^n p_{m,n}(x)$. 
By (1) we get that
{\footnotesize
\[
\alpha'(x) = \sum_{m=0}^n n \big( p_{m-1,n-1}(x) -
p_{m,n-1}(x) \big) = 
n \left( \sum_{m=1}^n p_{m-1,n-1}(x) - 
\sum_{m=0} p_{m,n-1}(x) \right) = 0. 
\] }
Since $\alpha(0)=1$ it follows that $\alpha(x) = 1$.

\item Put $\beta(x) = \sum_{m=0}^n m p_{m,n}(x)$.
By (1) and (2) we get that
{\footnotesize
\[
\beta'(x) = \sum_{m=0}^n mn \big( 
p_{m-1,n-1}(x) - p_{m,n-1}(x) \big)
\]
\[
= n \sum_{m=0}^{n-1} (m+1)p_{m,n-1}(x) - 
n \sum_{m=0}^{n-1} m p_{m,n-1}(x)
= n \sum_{m=0}^{n-1} p_{m,n-1}(x) = n.
\]
}
Since $\beta(0)=0$ it follows that $\beta(x) = nx$.

\item Put $\gamma(x) = \sum_{m=0}^n (nx-m)^2 p_{m,n}(x)$.
By (1), (2) and (3) we get that
{\footnotesize
\[
\gamma'(x) = \sum_{m=0}^n 2n (nx-m)p_{m,n}(x) +
n (nx-m)^2 \big( p_{m-1,n-1}(x) - p_{m,n-1}(x) \big)
\]
\[
= 2n^2 x \sum_m p_{m,n}(x) - 2n \sum_m m p_{m,n}(x)
+ n^3 x^2 \sum_m p_{m-1,n-1}(x) 
- n^3 x^2 \sum_m p_{m,n-1}(x)
\]
\[
- 2n^2 x \sum_m m ( p_{m-1,n-1}(x) - p_{m,n-1}(x) ) 
+ n \sum_m m^2 p_{m-1,n-1}(x) 
- n \sum_m m^2 p_{m,n-1}(x)
\]
\[
= 2n^2 x - 2n^2 x + n^3 x^2 - n^3 x^2
- 2n^2 x
+ n \sum_m \big( (m+1)^2-m^2 \big) p_{m,n-1}(x)
\]
\[
= -2n^2 x + 2n \sum_m m p_{m,n-1} + n \sum_m p_{m,n-1}
= -2n^2 x + 2n(n-1)x + n 
= n(1-2x).
\]
}
Since $\gamma(0)=0$ it follows that 
$\gamma(x) = nx(1-x)$.

\item The equality $F_n(0) = 0$ is clear. By (1) we get that:
{\footnotesize
\[
F_n'(x) = \sum_{m=0}^n \frac{f(m/n)}{n+1} 
\sum_{j=m+1}^{n+1} (n+1) \big( p_{j-1,n}(x) - p_{j,n}(x) 
\big) = \sum_{m=0}^n f \left( m/n \right) p_{m,n}(x).
\]
}

\item We use the argument in 
\cite[p. 8]{burkill1959}.
Take $\epsilon > 0$. By uniform continuity of $f$ 
it follows that there is $\delta > 0$ such that
$| f(x_1)-f(x_2) | < \epsilon/2$ whenever 
$x_1 , x_2 \in I$ and $|x_1 - x_2| < \delta$. 
Choose $n \geq 1$ so that 
$n > 4 || f ||/(\epsilon \delta^2)$. Take $x \in I$. 
The equalities (2) and (4) now yield that
{\footnotesize
\[
| f(x) - f_n(x) | \leq 
\sum_{|x-m/n| < \delta} | f(m/n)-f(x) | p_{m,n}(x) +
\sum_{|x-m/n| \geq \delta} | f(m/n)-f(x) | p_{m,n}(x)
\]
\[
\leq \frac{\epsilon}{2} \sum_m p_{m,n}(x) + 
\frac{2 ||f|| }{n^2 \delta^2} \sum_m (nx-m)^2 p_{m,n}(x) 
= \frac{\epsilon}{2} + \frac{2 ||f|| x(1-x)}{n \delta^2}
< \frac{\epsilon}{2} + \frac{\epsilon}{2} = \epsilon.
\qed
\]
}
\end{enumerate}


\noindent
\emph{Proof of theorem.}
The second part of our argument 
follows the proof of \cite[Thm. 9.13]{apostol1974}.
Let $x \in J$.
Since $F_m(0)=F_n(0)=0$ the mean value theorem
implies that there is $x_1$ between 
$0$ and $x$ which satisfies
$F_m(x) - F_n(x) = x (F_m'(x_1)-F_n'(x_1)) = 
x (f_m(x_1) - f_n(x_1)).$
Therefore $|| F_m - F_n || \leq ||f_m - f_n ||$
and thus 
$\{ F_n \}_{n=1}^\infty$ is uniformly convergent on $I$. 
By the proposition, the limiting function $F$
is continuous function on $I$.
Take $c \in J$. 
We wish to show that $F'(c)$ exists and is equal to $f(c)$.
To this end, define the polynomial $Q_n$ 
as the result of the polynomial division of 
$F_n(x) - F_n(c)$ by $x-c$. 
Note that $Q_n(c) = F_n'(c) = f_n(c)$.
Take $x \in I$ with $x \neq c$.
By the mean value theorem again, we get that 
\begin{center}
$Q_m(x) - Q_n(x) = 
\frac{F_m(x)-F_m(x) - (F_m(c)-F_n(c))}{x-c}
= f_m(x_2) - f_n(x_2)$ 
\end{center} 
for some $x_2$ between $x$ and $c$. In particular,
$|| Q_m - Q_n || \leq || f_m - f_n ||$
which shows that $\{ Q_n \}_{n=1}^\infty$
is uniformly convergent on $I$. 
By the proposition again the limiting function
$Q$ is continuous on $I$.
By letting $n \to \infty$ in the equality
$F_n(x) - F_n(c) = Q_n(x)(x - c)$ we get that
$F(x) - F(c) = Q(x)(x-c)$.
Continuity of $Q$ at $x = c$ implies that
$F'(c)$ exists and equals
$Q(c) = \lim_{n \to \infty} Q_n(c) = 
\lim_{n \to \infty} f_n(c) = f(c)$.
\qed

\end{document}